\documentclass[12pt]{article}
\usepackage{a4wide}
\usepackage{amssymb}
\usepackage{amsmath}
\usepackage[russian]{babel}

  \chardef\No=242
\newtheorem{theorem}{Теорема}{}
\newtheorem{definition}{Definition}{}
\newtheorem{corollary}{Corollary}{}
\newtheorem{lemma}{Lemma}{}
\newtheorem{remark}{Remark}{}

\begin{document}
\begin{center}
{\bf  Application of multidimensional Hardy operator and \\
its connection with a certain nonlinear differential equation\\
in weighted variable Lebesgue spaces }
\end{center}

\begin{center}
Rovshan A. Bandaliev \\
({\it Institute of Mathematics and Mechanics of National Academy of Sciences of Azerbaijan})\\
E-mail address: bandaliyevr@gmail.com
\end{center}
\vspace{2mm}

{\small  {\it Abstract.} In this paper a two weight criterion for multidimensional geometric mean operator
in variable exponent Lebesgue space is proved. Also, we found a criterion on weight functions expressing
one-dimensional Hardy inequality via a certain nonlinear differential equation. In particular, considered
nonlinear differential equation is nonlinear integro-differential equation.}

\vspace{2mm}

\vspace{2mm}

\vspace{2mm}
\begin{center}
{\bf  Introduction.}
\end{center}

It is well known that the variable exponent Lebesgue space appeared
in the literature for the first time already in [36].  Further development of this
theory was connected with the theory of modular function spaces. The
first systematic study of modular spaces is presented in [33]. In
the appendix, Nakano mentions explicitly variable exponent Lebesgue
spaces as an example of more general spaces he considers.
Somewhat later, a more explicit version of these spaces, namely
modular function spaces, were investigated by many mathematicians
(see [32]). The next step in the investigation of variable exponent spaces was
given in [42] and in [25]. The study of these spaces has been
stimulated by problems of elasticity, fluid dynamics, calculus of
variations and differential equations with non-standard growth
conditions$\,$ (see [39], [48]-[50]).

Inequalities are one of the most important instruments in many branches
of mathematics such as functional analysis, theory of differential and integral
equations, interpolation theory, harmonic analysis, probability theory  etc.
They are also useful in mechanics, physics and other sciences.
It is well known that the classical two weight inequality for the geometric
mean operator is closely connected to the one-dimensional Hardy inequality (see [15]).
Analogously, the P\'{o}lya-Knopp type inequalities with multidimensional geometric
mean operator is connected with multidimensional Hardy type operator. In the papers [6],
[44], [45] and [14] the connection of the Hardy inequality with a nonlinear differential
equation having a solution with certain special properties was considered. Therefore,
the consideration of this problems in variable exponent Lebesgue space is actual.

\begin{center}
\bf {1.$\;$ Preliminaries}
\end{center}

Let $R^{n}$ be the $n$-dimensional Euclidean space of points $x=\left(x_{1},..., x_{n}\right)$
and $\Omega$ be a Lebesgue measurable subset in $R^n$ and $\displaystyle{|x|= \left(\sum\limits_{i= 1}^n
x_i^2\right)^{1/2}}.$ Suppose that $p$ is a Lebesgue measurable function on $\Omega$ such that
$0< \underline p\le p(x)< \infty,$ $\underline p= ess\,\inf_{x\in \Omega} p(x)$  and $\omega$ is a weight function on $\Omega,$ i.e. $\omega$ is a non-negative, almost everywhere (a.e.) positive function on $\Omega.$ The Lebesgue measure of a set $\Omega$ will be denoted by $|\Omega|.$ It is well known that $\displaystyle{|B(0,1)|= \frac {\pi^{\frac n2}}{\Gamma\left(\frac n2+ 1\right)}},$ where $B(0,1)= \left\{x\in R^n;\;|x|< 1\right\}.$ Further, in this paper all sets and functions are supposed to be Lebesgue measurable. By $AC(0,\,\infty)$ we denote the set of absolutely continuous functions on $(0,\,\infty).$ For the sake of simplicity, the letter $C$ always denotes a positive constant which may change from one step to the next.

\begin{definition}
By $L_{p(x),\,\omega}(\Omega)$ we denote the set of measurable
functions $f$ on $\Omega$ such that for some $\lambda_0> 0$
$$
\int\limits_{\Omega} \left(\frac{|f(x)|\,\omega(x)}{\lambda_0}\right)^{p(x)}\,dx< \infty.
$$
Note that the expression
$$
\|f\|_{L_{p(x),\,\omega}(\Omega)}= \|f\|_{L_{p(\cdot),\,\omega}(\Omega)}=
\inf\left\{\lambda> 0:\;\;\int\limits_{\Omega} \left|\frac{f(x)\,\omega(x)}
{\lambda}\right|^{p(x)} \,dx\le 1\right\}.
$$
defines a quasi-Banach spaces. In particular, for $1\le p(x)< \infty$ the space $L_{p(x),\,\omega}(\Omega)$
is a Banach function space (see [11]) with respect to the expression $\|f\|_{L_{p(x),\,\omega}(\Omega)}.$
\end{definition}

For $\omega= 1$ the space $L_{p(x),\omega}(\Omega)$ coincides with the variable Lebesgue space $L_{p(x)}(\Omega).$

We reduce two examples which characterize the norm of this space. \\
{\bf Example 1.} {\it Let $\;p(x)= \left\{
            \begin{array}{l}
            2 \quad\quad for\quad x\in  \Omega\\
            3 \quad\quad for\quad x\in R^n\setminus \Omega\\
            \end{array}
            \right.$ and $f\in L_2\left(R^n\right)\bigcap L_3\left(R^n\right).$}

We calculate $\|f\|_{L_{p(x)}(R^n)}.$ By the definition we have
$$
\|f\|_{L_{p(x)}(R^n)}= \inf\left\{\lambda> 0:\;\;\int\limits_{\Omega}
\left|\frac{f(x)}{\lambda}\right|^2 \,dx+ \int\limits_{R^n\setminus\Omega}
\left|\frac{f(x)}{\lambda}\right|^3 \,dx \le 1\right\}=
$$
$$
= \inf\left\{\lambda> 0:\;\;\frac{a_1}{\lambda^2}+ \frac{a_2}{\lambda^3}\le 1\right\}=
\inf\left\{\lambda> 0:\;\;\lambda^3- a_1\,\lambda- a_2\ge 0\right\},
$$
where $\displaystyle{a_1= \int\limits_{\Omega} f^2(x)\;dx}$ and
$\displaystyle{a_2= \int\limits_{R^n\setminus\Omega} |f(x)|^3\;dx}.$
Now we solve the inequality $\;\lambda^3- a_1\,\lambda- a_2\ge 0.$ We
consider three different  cases.\\
{\bf Case 1.} Let $\displaystyle{\frac{a_2^2}4- \frac{a_1^3}{27}> 0}.$ Then the cubic equation
$\,\lambda^3- a_1\,\lambda- a_2= 0$ has one real root and two complex conjugate roots.
Namely, $\displaystyle{\lambda_1= \sqrt[3]{\frac {a_2}2+ \sqrt{\frac{a_2^2}4- \frac{a_1^3}{27}}}+
\sqrt[3]{\frac {a_2}2- \sqrt{\frac{a_2^2}4- \frac{a_1^3}{27}}}},$ $\displaystyle{\lambda_2= -\frac{\lambda_1}2+
i\,\sqrt{3}\,\frac{\sqrt[3]{\frac {a_2}2+ \sqrt{\frac{a_2^2}4- \frac{a_1^3}{27}}}-
\sqrt[3]{\frac {a_2}2- \sqrt{\frac{a_2^2}4- \frac{a_1^3}{27}}}}2}$ and $\lambda_3= \overline \lambda_2.$ It is obvious that
$\lambda^3- a_1\,\lambda- a_2= \left(\lambda- \lambda_1\right)\left(\lambda^2+ \lambda_1\,\lambda+ \left|\lambda_2\right|^2\right)$
and $\lambda^2+ \lambda_1\,\lambda+ \left|\lambda_2\right|^2> 0$ for all $\lambda\in (-\infty, +\infty).$ Therefore the inequality
$\;\lambda^3- a_1\,\lambda- a_2\ge 0$ holds if and only if $\lambda\ge \lambda_1$ and
$$
\|f\|_{L_{p(x)}(R^n)}= \sqrt[3]{\frac {a_2}2+ \sqrt{\frac{a_2^2}4- \frac{a_1^3}{27}}}+
\sqrt[3]{\frac {a_2}2- \sqrt{\frac{a_2^2}4- \frac{a_1^3}{27}}}.
$$
\\
{\bf Case 2.} Let $\displaystyle{\frac{a_2^2}4- \frac{a_1^3}{27}= 0}.$ Then $\displaystyle{\lambda^3- a_1\,\lambda- a_2=
\left(\lambda- 2\,\sqrt[3]{\frac{a_2}2}\right)\left(\lambda+ \sqrt[3]{\frac{a_2}2}\right)^2}$ and the inequality
$\;\lambda^3- a_1\,\lambda- a_2\ge 0$ holds if and only if $\displaystyle{\lambda\ge 2\,\sqrt[3]{\frac{a_2}2}}$ and
$\displaystyle{\|f\|_{L_{p(x)}(R^n)}= 2\,\sqrt[3]{\frac{a_2}2}}.$
\\
{\bf Case 3.} Let $\displaystyle{\frac{a_2^2}4- \frac{a_1^3}{27}< 0}.$ Then the equation $\displaystyle{\lambda^3-
a_1\,\lambda- a_2= 0}$ has three distinct real roots. We denote by $\alpha_1,$ $\alpha_2$ and $\alpha_3$ the roots
of this equation. By Vi\`{e}te's formulas one root of this equation is positive and two roots are negative. Let $\alpha_1> 0.$
Then $\displaystyle{\lambda^3- a_1\,\lambda- a_2= }$\linebreak $\displaystyle{\left(\lambda- \alpha_1\right)\left(
\lambda^2+ \alpha_1\,\lambda+ \alpha_1^2- a_1\right)= 0},$ $\displaystyle{\alpha_2= \frac{-\alpha_1+ \sqrt{4\,a_1- 3\,
\alpha_1^2}}2}$ and $\displaystyle{\alpha_3= \frac{-\alpha_1- \sqrt{4\,a_1- 3\,\alpha_1^2}}2},$ where $\displaystyle{
\sqrt{a_1}< \alpha_1< \frac 2{\sqrt 3}\,\sqrt{a_1}}.$ It is obvious that $\alpha_3< \alpha_2< \alpha_1.$ Therefore
the inequality $\lambda^3- a_1\,\lambda- a_2\ge 0$ holds if and only if $\lambda\in \left[b_3,\;b_2\right]\bigcup
\left[\alpha_1, \infty\right)$ and by the definition of the norm we have $\lambda\ge \alpha_1.$ Thus,
$\displaystyle{\|f\|_{L_{p(x)}(R^n)}= \alpha_1}.$
\\
{\bf Example 2.} {\it Let $n= 1,$ $x\in [1, \,\infty),$ $\;p(x)= x$ and $f(x)= 1.$}

We calculate $\|1\|_{L_{p(x)}([1,\infty))}.$ We have
$$
\|1\|_{L_{p(x)}([1,\infty))}= \inf\left\{\lambda> 0:\;\;\int\limits_{1}^\infty
\frac 1{\lambda^x}\,dx \le 1\right\}.
$$
It is obvious that $\displaystyle{\int\limits_{1}^\infty \frac 1{\lambda^x}\,dx= \frac{1}{\lambda\,\ln \lambda}}\,$ if
$\;\lambda> 1\;$ and
$$
\inf\left\{\lambda> 0:\;\;\int\limits_{1}^\infty \frac 1{\lambda^x}\,dx \le 1\right\}= \inf\left\{\lambda> 1:\;\;
\frac{1}{\lambda\,\ln \lambda}\le 1\right\}= \inf\left\{\lambda> 1:\;\;
\lambda^{\lambda}\ge e\right\}.
$$
Thus, $\|1\|_{L_{p(x)}([1,\infty))}= 1,7712...$

The following theorems are known.
\begin{theorem}{[3]} Let $1\le \underline p\le p(x)\le q(y)\le \overline q<
\infty$ for all $x\in \Omega_1\subset R^n$ and $y\in \Omega_2\subset
R^m.$  If $p(x)\in C\left(\Omega_1\right),$ then the inequality
$$
\left\|\|f\|_{L_{p(\cdot)}\left(\Omega_1\right)}\right\|_{L_{q(\cdot)}\left(\Omega_2\right)}\le
\left(\frac{\overline p}{\underline q}+ \frac{\overline q- \underline p}{\overline q}\right)^{\frac 2{\underline p}}
\,\left\|\|f\|_{L_{q(\cdot)}\left(\Omega_2\right)}\right\|_{L_{p(\cdot)}\left(\Omega_1\right)}
$$
is valid, where $\underline q= \mbox {ess}\,\inf\limits_{\Omega_2}
q(x),$ $\overline q= \mbox {ess}\,\sup\limits_{\Omega_2} q(x)$ and
$C\left(\Omega_1\right)$ is the space of continuous functions in
$\Omega_1$ and $f:\Omega_1\times \Omega_2\rightarrow R$ is any
measurable function such that
$$
\left\|\|f\|_{q,\Omega_2}\right\|_{p,\Omega_1}= \inf \left\{ \mu>
0:\;\;\int\limits_{\Omega_1}\left(\frac{\|f(x,\cdot)\|_{q(\cdot),\Omega_2}}
{\mu}\right)^{p(x)}\,dx\le 1\right\}< \infty.
$$
\end{theorem}

Let $\displaystyle{Hf(x)= \int\limits_{|y|< |x|} f(y)\,dy},$ where $f\ge 0$ and $B(0, |x|)= \left\{y\in R^n;\; |y|< |x|\right\}.$
Now we formulate the criteria on boundedness of multidimensional Hardy type operator in weighted variable Lebesgue spaces.
\begin{theorem}{[5]} Let $q(x)$ be a measurable function on $R^n,$ $1< p\le q(x)\le \overline q< \infty$ and $\displaystyle{p'=\frac{p}{p- 1}}.$ Suppose that $v(x)$ and $w(x)$ are weights on $R^n.$ Then the inequality
$$
\left\|Hf\right\|_{L_{q(\cdot),\, w}\left(R^n\right)}\le C
\,\left\|f\right\|_{L_{p,\,v}(R^n)} \eqno(1.1)
$$
holds, for every $f\ge 0$ if and only if there exists $\alpha\in (0, 1)$ such that
$$
A(\alpha, p, q)= \sup\limits_{t> 0}\left(\int\limits_{|y|< t} v^{-p'}(y)\;dy\right)^
{\frac{\alpha}{p{\,'}}}\left\|\left(\int\limits_{|y|< |\cdot|} v^{-p'}(y)\;dy\right)^
{\frac{1- \alpha}{p^{\,'}}}\right\|_{L_{q(\cdot),\, w}(|x|> t)}< \infty. \eqno(1.2)
$$
Moreover, if $C> 0$ is the best possible constant in (1.1), then
$$
\sup\limits_{0< \alpha< 1} \frac{p{\,'}\,A(\alpha, p,
q)}{(1- \alpha)\,\left[\left(\frac{p{\,'}}{1-\alpha}\right)^{p}+ \frac
1{\alpha\left(p- 1\right)}\right]^{1/ p}} \le
C\le
$$
$$
\le \left(\frac{p}{\underline q}+ \frac{\overline q- p}{\overline q}\right)^{\frac 2{
p}}\,\inf\limits_{0< \alpha< 1} \frac{A(\alpha, p, q)}{(1- \alpha)^{1/p{\,'}}}.
$$
\end{theorem}

\begin{remark} Note that Theorem 3 in the case $n= 1,$  $q(x)= q= const$ for $x\in (0, \infty)$ and $\alpha=
\frac{s- 1}{p- 1}$ ($1< s< p-1$) was proved in [46] and in multidimensional Hardy type operators it was
proved in [5] (see also [4]). Mainly the multidimensional Hardy operator it was investigated in the papers
[2], [8], [9], [12], [13], [18], [19], [26], [29], [31], [34], [40], [41], [43], [47]  etc. The prehistory of the
Hardy inequality with different $p$ and $q$ were started in [7] and independently derived in [23], [30] and [1].
Two-weighted criterion for one-dimensional Hardy operator in weighted variable $L_{p(x),\,w}([0, 1])$ spaces was proved in [24]. The criteria on the exponent function for validity the boundedness of one-dimensional Hardy operator in variable $L_{p(x)}(0,\, \infty)$ spaces was proved in [16]. Also, other type two-weight criteria for the multidimensional Hardy operator in weighted variable Lebesgue spaces was proved in [28].
\end{remark}

\vspace{5mm}
\begin{center}
{\bf 2. Main results.}
\end{center}

We consider the multidimensional geometric mean operator defined as
$$
Gf(x)= \exp\left(\frac{1}{|B(0,\,|x|)|}\int\limits_{B(0,\,|x|)} \ln\,f(y)\,dy\right),
$$
where $f> 0$ and $|B(0,\,|x|)|= |B(0,1)|\,|x|^n.$ It is obvious that $G\left(f_1\cdot f_2\right)(x)= $
\linebreak $G f_1(x)\cdot Gf_2(x).$

Now we formulate a two-weight criterion on boundedness of multidimensional geometric mean operator in variable Lebesgue spaces.

\begin{theorem} Let $q(x)$ be a measurable function on $R^n$ and $0< p\le q(x)\le \overline q< \infty.$ Suppose that
$v(x)$ and $w(x)$ are weights on $R^n.$ Then the inequality
$$
\left\|Gf\right\|_{L_{q(\cdot),\, w}\left(R^n\right)}\le C
\,\left\|f\right\|_{L_{p,\,v}(R^n)} \eqno(2.1)
$$
holds, for every $f> 0$ if and only if there exists $s\in (1, p)$ such that
$$
D(s, p, q)= \!\sup\limits_{t> 0} |B(0, t)|^{\frac{s- 1} p}\left\|\frac{w(\cdot)}{|B(0, |\cdot|)|^{\frac {s}p}}\!\!
\,\exp\left(\frac 1{|B(0, |\cdot|)|}\int\limits_{B(0,|\cdot|)} \!\!\ln\frac 1{v(y)}\,dy\right)\right\|_
{L_{q(\cdot)}(|x|> t)}\!\!\!\!\!\!\!\!\!\!\!\!\!\!\!\!\!\!\!< \infty. \eqno(2.2)
$$
Moreover, if $C> 0$ is the best possible constant in (2.1), then
$$
\sup\limits_{s> 1} \frac{e^{\frac sp}}{\left(e^s+ \frac 1{s- 1}\right)^{1/p}}\,
D(s, p,q) \le C \le \left(\frac{p}{\underline q}+ \frac{\overline q - p}{\overline q}\right)^{\frac {2}{p}}\,
\inf\limits_{s> 1} e^{\frac{s- 1}{p}}\, D(s, p, q).
$$
\end{theorem}

{\bf Proof.} Let $\displaystyle{\alpha= \frac{s- 1}{p- 1}},$ where $1< s< p.$  We replace  $f$ with $f^{\beta},$ $v$ with $v^{\beta},$ $w$ with $\displaystyle{\frac {w^{\beta}(x)}{|B(0, |x|)|}},$ $0< \beta< p,$ and  $\,p\,$ with $\displaystyle{\frac p{\beta}}$ and $\,q(x)\,$ with $\displaystyle{\frac {q(x)}{\beta}}$ in (1.1), (1.2). We find that for $\displaystyle{1< s< \frac p{\beta}}$
$$
\left\|\frac{w^{\beta}}{|B(0, |\cdot|)|}\,H(f^{\beta})\right\|_{L_{\frac{q(\cdot)}{\beta}}(R^n)}= \left\|\left(\frac 1{|B(0, |\cdot|)|}
\int\limits_{B(0,\,|\cdot|)} f^{\beta}(y)\,dy\right)^{1/\beta}\right\|^{\beta}_{L_{q(\cdot),\, w}\left(R^n\right)}\le
$$
$$
\le C_{\beta}\left(\int\limits_{R^n} [f(y) v(y)]^{p}\,dy\right)^{\beta/p}.
$$
Then the inequality
$$
\left\|\left(\frac 1{|B(0, |\cdot|)|} \int\limits_{B(0,\,|\cdot|)} f^{\beta}(y)\,dy\right)^{1/\beta}
\right\|_{L_{q(\cdot),\, w}\left(R^n\right)}\le C_{\beta}^{1/\beta}\left(\int\limits_{R^n} [f(y) v(y)]^{p}\,dy\right)^{1/p} \eqno(2.3)
$$
holds if and only if
$$
A\left(\frac{s- 1}{p- 1}, \frac p{\beta}, \frac q{\beta}\right)=
$$
$$
=\left[\sup\limits_{t> 0} \left(\int\limits_{|y|< t} [v(y)]^{-\frac{\beta\,p}{p- \beta}}\,dy\right)^
{\frac{s- 1}p}\left\|\left(\frac 1{|B(0, |\cdot|)|^{\frac{p}{p- \beta s}}}
\int\limits_{|y|< |\cdot|} [v(y)]^{-\frac{\beta p}{p- \beta}}\,dy\right)^
{\frac{p- \beta s}{\beta p}} \right\|_{L_{q(\cdot),\, w}(|x|> t)}\right]^{\beta}=
$$
$$
= B^{\beta}\left(s, p,  q, \beta\right)< \infty
$$
and
$$
\sup\limits_{1< s< \frac p\beta} \left[\frac{\left(\frac{p}{p- s\beta}\right)^{\frac p\beta}}
{\left(\frac{p}{p- s\beta}\right)^{\frac p\beta}+ \frac 1{s- 1}}\right]^{\beta/p} B^{\beta}\left(s, p,  q,
\beta\right)\le C_{\beta}\le
$$
$$
\le \left(\frac{p}{\underline q}+ \frac{\overline q - p}{\overline q}\right)^{\frac {2\,\beta}
{p}}\,\inf\limits_{1< s< \frac p\beta} \left(\frac{p- \beta}{p- s\beta}\right)^{\frac {p- \beta}{p}}
\,B^{\beta}\left(s, p,  q, \beta\right), \eqno(2.4)
$$
where $\underline q$ is replaced by $\displaystyle{\frac {\underline q}{\beta}}$ and $\overline q$ is replaced by
$\displaystyle{\frac {\overline q}{\beta}}.$
By L'Hospital rule, we get
$$
\lim\limits_{\beta\to +0} \left(\frac 1{|B(0, |x|)|^{\frac{p}{p- \beta s}}}\int\limits_{|y|< |x|} [v(y)]^{-\frac{\beta p}
{p- \beta}}\,dy\right)^{\frac{p- \beta s}{\beta p}}=
$$
$$
= \lim\limits_{\beta\to +0} \exp\left[\frac{p\,\ln\frac 1{|B(0, |x|)|}
+ (p- \beta s)\,\ln\left(\int\limits_{|y|< |x|}[v(y)]^{-\frac{\beta p}{p- \beta}}\,dy\right)}{p\,\beta}\right]=
$$
$$
= \lim\limits_{\beta\to +0} \exp\left[-\frac sp\,\ln\left(\int\limits_{|y|< |x|}[v(y)]^{-\frac{\beta p}{p- \beta}}\,
dy\right)+ \frac{(p- \beta s)\left(\frac{p}{p- \beta}\right)^2\int\limits_{|y|< |x|} [v(y)]^{-\frac{\beta p}{p- \beta}}
\ln \frac 1{v(y)}\,dy}{p\,\int\limits_{|y|< |x|}[v(y)]^{-\frac{\beta p}{p- \beta}}\,dy}\right]=
$$
$$
= \exp\left[\frac sp\,\ln\frac 1{|B(0, |x|)|}+ \frac{\int\limits_{|y|< |x|} \ln \frac 1{v(y)}\,dy}{|B(0, |x|)|}\right]=
\frac 1{|B(0, |x|)|^{\frac sp}}\,\exp\left(\frac 1{|B(0, |x|)|}\int\limits_{B(0,|x|)} \ln\frac 1{v(y)}\,dy\right).
$$
Therefore
$$
\lim\limits_{\beta\to +0} B\left(s, p, q, \beta\right)=
$$
$$
= \sup\limits_{t> 0} |B(0, t)|^{\frac{s- 1} p} \left\|\frac{w(\cdot)}{|B(0, |\cdot|)|^{\frac {s}p}}
\,\exp\left(\frac 1{|B(0, |\cdot|)|}\int\limits_{B(0,|\cdot|)} \ln\frac 1{v(y)}\,dy\right)\right\|_
{L_{q(\cdot)}(|x|> t)}=
$$
$$
= D(s, p, q)< \infty
$$
and
$$
\sup\limits_{s> 1} \frac{e^{\frac sp}}{\left(e^s+ \frac 1{s- 1}\right)^{1/p}}\,D(s, p,q)\le
\lim\limits_{\beta\to +0} C_{\beta}^{1/\beta}\le \left(\frac{p}{\underline q}+ \frac{\overline q - p}
{\overline q}\right)^{\frac {2}{p}}\,\inf\limits_{s> 1}e^{\frac {s- 1}{p}}\,D(s, p, q). \eqno(2.5)
$$
Further, we have
$$
\lim\limits_{\beta\to +0} \left(\frac 1{|B(0, |x|)|}\int\limits_{B(0,\,|x|)} f^{\beta}(y)\,dy\right)^{1/\beta}=
\exp\left(\frac{1}{|B(0,\,|x|)|}\int\limits_{B(0,\,|x|)} \ln\,f(y)\,dy\right)= Gf(x).
$$
Formula (2.4) implies $\displaystyle{\lim\limits_{\beta\to +0} C_{\beta}= 1},$ and according to (2.2) and (2.5)
$\displaystyle{\lim\limits_{\beta\to +0} C_{\beta}^{1/\beta}= C< \infty}.$ Therefore the inequality (2.3) is valid.
Moreover, from (2.3) for $\beta\to +0$ we obtain
$$
\left\|Gf\right\|_{L_{q(\cdot),\, w}\left(R^n\right)}\le C
\,\left\|f\right\|_{L_{p,\,v}(R^n)}
$$
and by (2.5)
$$
\sup\limits_{s> 1} \frac{e^{\frac sp}}{\left(e^s+ \frac 1{s- 1}\right)^{1/p}}\,
D(s, p,q) \le C \le \left(\frac{p}{\underline q}+ \frac{\overline q - p}{\overline q}\right)^{\frac {2}{p}}\,
\inf\limits_{s> 1} e^{\frac{s- 1}{p}}\, D(s, p, q).
$$

This completes the proof of Theorem 3.
\begin{remark}
Let $q(x)= q= const$ and $n= 1.$ Note that the simplest case of (2.1) with $v= w= 1$ and $p= q= 1$
was considered in [15], and in [22]. Later on, this inequality was generalized in various ways by
many authors in [10], [17], [20], [21], [27], [35], [37], [38], [46]  etc.
\end{remark}

\begin{corollary}
Let $q(x)= q= const,$ $0< p\le q< \infty$ and let $f$ be a positive function on $R^n.$ Then
$$
\left(\int\limits_{R^n} [Gf(x)]^q\,|x|^{\gamma\,q}\;dx\right)^{1/q}\le
C\;\left(\int\limits_{R^n} f^p (x)\,|x|^{\beta\,p}\;dx\right)^{1/p} \eqno(2.6)
$$
holds with a finite constant $C$ if and only if
$$
\gamma+ \frac nq= \frac{\beta}n+ \frac np
$$
and the best constant $C$ has the following condition:
$$
\sqrt[q]{\frac p{nq}}\;e^{\frac{\beta}{n^2}}\;|B(0,\,1)|^{\frac 1q- \frac 1p}\,\sup\limits_{s> 1}
\frac{e^{\frac sp}\,(s- 1)^{\frac 1p- \frac 1q}}{\left[(s- 1)e^s+ 1\right]^{1/p}}\le C\le
\frac{|B(0,\,1)|^{\frac 1q- \frac 1p}\;e^{\frac{\beta}{n^2}+ \frac 1q}}{\sqrt [q]{n}}.
$$
\end{corollary}

\begin{remark}
Note that if $p= q,$ then the inequality (2.6) is sharp with the constant $\displaystyle{C= \frac {e^{\frac{\beta}{n^2}+ \frac 1p}}
{\sqrt [p]{n}}}.$
\end{remark}

\begin{corollary} Let $x\in R^n,$ $0< p\le q(x)< \infty,\,$ $q(x)= \left\{\begin{array}{l}
1 \quad\, for\quad |x|< 1\\
2\quad for\quad |x|\ge 1,\\
\end{array}
\right.$ and let $f$ be a positive function on $R^n.$ Suppose that $v(x)= 1$ and $w(x)= |x|^{\beta}.$ Then
$$
\left\|Gf\right\|_{L_{q(\cdot),\, |\cdot|^{\beta}}\left(R^n\right)}\le
C\;\left(\int\limits_{R^n} f^p(x)\,dx\right)^{1/p}
$$
holds with a finite constant $C$ if and only if
$$
n\left(\frac sp- 1\right)\le \beta\le n\left(\frac 1p- \frac 12\right), \;\;s\in \left(1, 1+ \frac p2\right]
$$
and the best constant $C$ has the following condition:
$$
\sup\limits_{1< s\le 1+ \frac p2} \frac{e^{\frac sp}}{\left(e^s+ \frac 1{s- 1}\right)^{1/p}}\,
D'(s, p,q) \le C \le \left(\frac{p}{\underline q}+ \frac{\overline q - p}{\overline q}\right)^{\frac {2}{p}}\,
\inf\limits_{1< s\le 1+ \frac p2} e^{\frac{s- 1}{p}}\, D'(s, p, q),
$$
where $\displaystyle{ D'(s, p, q)= |B(0, 1)|^{-\frac 1p}\, \sup\limits_{t> 0} t^{\frac{n(s- 1)} p}
\,\left\|\,|\cdot|^{\beta- \frac {ns}p}\,\right\|_{L_{q(\cdot)}(|\cdot|> t)}< \infty.}$
\end{corollary}

Now we consider an application in the theory of nonlinear ordinary differential equation. Let $L(t, \omega, y)=
\left\|\omega\,y^{1/p'}\right\|_{L_{q(x)}(x> t)},$ where $t\in (0, \infty)$ and $\omega$ is a
weight function defined on $(0,\infty).$

\begin{lemma} Let $1< p\le q(x)\le \overline q< \infty.$ Suppose that $\omega_1$
and $\omega_2$ are weight functions defined on $(0, \infty).$ Let the equation
$$
L\left(t,\omega_2,y\right)- \lambda\omega_1(t) \left(y'(t)\right)^{1/p'}= 0, \;\;(\lambda> 0) \eqno(2.7)
$$
have a solution $\;y\;$ such that
$$
y(t)> 0,\,\;y'(t)> 0,\;\, y\in AC(0,\,\infty). \eqno(2.8)
$$

Then the weighted norm inequality
$$
\|u\|_{L_{q(\cdot), \omega_2}(0, \infty)}\le \lambda\,\left(\frac{p}{\underline q}+ \frac{\overline q- p}
{\overline q}\right)^{\frac 2{p}}\,\|u'\|_{L_{p, \omega_1}(0, \infty)}
$$
holds, where $u\in AC(0,\,\infty)$ and $u(0)= \lim\limits_{t\to +0} u(t)= 0.$
\end{lemma}

{\bf Proof.} Applying the H\"{o}lder inequality we have
$$
u(x)= \int\limits_0^x u'(t)\,dt= \int\limits_0^x u'(t)(y'(t))^{-\frac 1{p'}}\,(y'(t))^{\frac 1{p'}}\,dt\le
$$
$$
\le \left(\int\limits_0^x y'(t)\,dt\right)^{\frac 1{p'}}\,\left\|u'\,\left[y'\right]^{- \frac 1{p'}}\right\|_{L_{p}(0,\, x)}\le
[y(x)]^{\frac 1{p'}}\,\left\|u'\,\left[y'\right]^{- \frac 1{p'}}\right\|_{L_{p}(0,\, x)}.
$$
Thus
$$
\|u\|_{L_{q(\cdot),\, \omega_2}(0, \infty)}\le \left\|\omega_2(\cdot)[y(\cdot)]^{\frac 1{p'}}\,\|u'\,\left[y'\right]^{- \frac 1{p'}}\|_
{L_{p}(0,\, \cdot)}\right\|_{L_{q(\cdot)}(0, \infty)}=
$$
$$
= \left\|\left\|\omega_2(\cdot)[y(\cdot)]^{\frac 1{p'}}\,u'\,\left[y'\right]^{- \frac 1{p'}}\,\chi_{(0,\;\cdot)}\right\|_{L_{p}(0, \infty)}
\right\|_{L_{q(\cdot)}(0, \infty)}.
$$

By using Theorem 1 we have
$$
\left\|\left\|\omega_2(\cdot)[y(\cdot)]^{\frac 1{p'}}\,u'\,\left[y'\right]^{- \frac 1{p'}}\,\chi_{(0,\;\cdot)}\right\|_{L_{p}(0, \infty)}
\right\|_{L_{q(\cdot)}(0, \infty)}\le
$$
$$
\le \left(\frac{p}{\underline q}+ \frac{\overline q- p}{\overline q}\right)^{\frac 2{p}}\,\left\|\left\|\omega_2(\cdot)[y(\cdot)]^
{\frac 1{p'}}\,u'\,\left[y'\right]^{- \frac 1{p'}}\,\chi_{(0,\;\cdot)}\right\|_{L_{q(\cdot)}(0, \infty)}\right\|_{L_{p}(0, \infty)}=
$$
$$
= \left(\frac{p}{\underline q}+ \frac{\overline q- p}{\overline q}\right)^{\frac 2{p}}\,\left\|\left\|\omega_2\,
y^{\frac 1{p'}}\right\|_{L_{q(\cdot)}(t, \infty)}\,u'\,\left[y'\right]^{- \frac 1{p'}}\right\|_{L_{p}(0, \infty)}=
$$
$$
= \lambda\,\left(\frac{p}{\underline q}+ \frac{\overline q- p}{\overline q}\right)^{\frac 2{p}}\,\left\| \omega_1\,
\left[y'\right]^{\frac 1{p'}}\,u'\,\left[y'\right]^{- \frac 1{p'}}\right\|_{L_{p}(0, \infty)}=
$$
$$
=\lambda\,\left(\frac{p}{\underline q}+ \frac{\overline q- p}{\overline q}\right)^{\frac 2{p}}\,\left\| \omega_1\,u'\right\|_
{L_{p}(0, \infty)}= \lambda\,\left(\frac{p}{\underline q}+ \frac{\overline q- p}{\overline q}\right)^{\frac 2{p}}\,\left\| u'\right\|_
{L_{p,\,\omega_1}(0, \infty)}.
$$
This proves the Lemma 1.

Put
$$
K=  p'\,\inf\,\sup\limits_{x> 0} \frac 1{f(x)- x- p'\int\limits_0^t \frac{\omega'_1(s)\,f(s)}{\omega_1(s)}\,ds} \int\limits_0^x
\frac{[f(t)]^{1/p'+ 1}  P\left(t, \omega_2, y, f\right)}{\omega_1(t)\,[y(t)]^{1/p'}}\,dt, \eqno(2.9)
$$
where $a$ is a fixed number in $(0, \infty),$ $P\left(t, \omega_2, y, f\right)\ge 0$ for all $t> 0,$
$y(t)$ is a fixed positive solution of equation (2.7) and the infimum is taken over the class of measurable functions such that
$$
f(x)> x+  p'\int\limits_0^t \frac{\omega'_1(s)\,f(s)}{\omega_1(s)}\,ds\quad \mbox{for all}\quad x> 0.
$$

The following lemma gives the relation between the number $K$ and the problem (2.7), (2.8).
\begin{lemma}
Let $\lambda> 0$ be the number from Lemma 1 and let $K$ be given by (2.9). Suppose that $\omega_1$
and $\omega_2$ are weight functions defined on $(0, \infty)$ and the derivative
$\omega'_1(t)$ exists for all $t\in (0, \infty)$ and $\omega_1(t)\ge \omega_1(0)> 0.$

Then the following statements are equivalent:

(i)$\,$ if the problem (2.7), (2.8) has a solution with a locally absolutely continuous first derivative, then $\lambda \ge K;$

(ii)$\,$ if $K< +\infty,$ then the problem (2.7), (2.8) has a solution for every $\lambda > K.$
\end{lemma}

{\bf Proof.} Assume that $(i)$ holds. Let $y_0(x)$ be a solution of (2.7), (2.8). Let us take $\displaystyle{w= \frac {y_0}{y'_0}}.$
Assume that $\displaystyle{P\left(t, \omega_2, y_0, w\right)= -\frac{d}{dt}L\left(t, \omega_2, y_0\right)= P(t)}.$ It is obvious
that $P(t)\ge 0$ for all $t> 0.$ Then by virtue of (2.7) $w$ is a positive solution of the equation
$$
w'(t)= \frac{p'\,\omega_1'(t)\,w(t)}{\omega_1(t)}+ \frac{p' [w(t)]^{1/p'+ 1}
P(t)}{\lambda\,\omega_1(t)\,\left(y_0(t)\right)^{1/p'}}+ 1. \eqno(2.10)
$$
Hence (2.10) implies
$$
w(t)\ge \int\limits_0^t w'(s)\,ds= p'\,\int\limits_0^t \frac{\omega'_1(s)\,w(s)}{\omega_1(s)}\,ds+ \frac{p'}{\lambda}
\int\limits_0^t \frac{[w(s)]^{1/p'+ 1}\,P(s)}{\omega_1(s)\,\left(y_0(s)\right)^{1/p'}}\,ds+ t. \eqno(2.11)
$$
From (2.11) implies that
$$
\lambda\ge \frac {p'}{w(t)- t- p'\int\limits_0^t \frac{\omega'_1(s)\,w(s)}{\omega_1(s)}\,ds} \int\limits_0^t
\frac{[w(s)]^{1/p'+ 1}\,P(s)}{\omega_1(s)\,\left(y_0(s)\right)^{1/p'}}\,ds. \eqno(2.12)
$$
From (2.11), (2.12) and (2.9) it follows that $\lambda\ge K$ and the proof of $(i)\Rightarrow (ii)$ is complete.

Assume now $(ii)$ holds. Let us fix $\lambda> K.$
By the definition of $K$ there exists a measurable function $f(x)$ such that
$$
f(x)\ge x+  p'\int\limits_0^x \frac{\omega'_1(t)\,f(t)}{\omega_1(t)}\,dt+ \frac{p'}{\lambda}\int\limits_0^x
\frac{[f(t)]^{1/p'+ 1} P(t)} {\omega_1(t)\,\left[y(t)\right]^{1/p'}}\,dt.\eqno(2.13)
$$
Let us define a sequence $w_n(x)$ by setting
$$
w_0(x)= f(x),
$$
$$
w_{n}(x)= x+ p'\int\limits_0^t \frac{\omega'_1(s)\,w_{n- 1}(s)}{\omega_1(s)}\,ds+ \frac{p'}{\lambda}\int\limits_0^x
\frac{[w_{n- 1}(t)]^{1/p'+ 1} P_{n- 1}(t)}{\omega_1(t)\,\left[y(t)\right]^{1/p'}}\,dt\;
(n= 1,2,\ldots), \eqno(2.14)
$$
where  $P_0(t)= P(t)$ and $P_n(t)\ge 0$ for all $n\in \Bbb {N}.$ From (2.13) it follows that $w_0(x)\ge w_1(x).$
We put $w_{n- 1}(x)\ge w_n(x)$ and let $P_{n}(t)$ be decreasing sequences with respect to $n$ on $(0, \infty),$
where $n\in \Bbb {N}.$ Then
$$
\int\limits_0^x \frac{[w_{n- 1}(t)]^{1/p'+ 1} P_{n- 1}(t)}{\omega_1(t)\,\left[y(t)\right]^{1/p'}}\,
dt\ge \int\limits_0^x \frac{[w_{n}(t)]^{1/p'+ 1} P_n(t)}{\omega_1(t)\,\left[y(t)\right]^{1/p'}}\,dt
$$
and
$$
w_{n}(x)- w_{n+ 1}(x)\ge p'\int\limits_0^x \frac{\omega'_1(s)\,\left[w_{n- 1}(s)- w_{n}(s)\right]}{\omega_1(s)}\,ds\ge
$$
$$
\ge \inf\limits_{s\in (0,\,\infty)} \left[w_{n- 1}(s)- w_{n}(s)\right]\,\int\limits_0^x \frac{\omega'_1(s)}{\omega_1(s)}\,ds
= \inf\limits_{t\in (0,\,\infty)} \left[w_{n- 1}(t)- w_{n}(t)\right]\,\ln \frac{\omega_1(x)}{\omega_1(0)}\,\ge 0.
$$
Since $w_n(x)\ge 0,$ the sequence (2.14) converges. We denote its limit by $w(x).$
By the Levi monotone convergence theorem it follows that $w$ is a nonnegative solution of the equation
$$
w(x)= x+ p'\int\limits_0^x \frac{\omega'_1(t)\,w(t)}{\omega_1(t)}\,dt+ \frac{p'}{\lambda}\int\limits_0^x
\frac{[w(t)]^{1/p'+ 1} P(t)}{\omega_1(t) \,\left[y(t)\right]^{1/p'}}\,dt,
$$
where $\displaystyle{P(t)= \lim\limits_{n\to \infty} P_n(t)}.$ Hence, $w$ is absolutely continuous and
satisfies the equation
$$
w'(x)= 1+ \frac{p'\,\omega'_1(x)\,w(x)}{\omega_1(x)}+ \frac{p'}{\lambda}\,
\frac{[w(x)]^{1/p'+ 1} P(x)}{\omega_1(x)\,\,\left[y(x)\right]^{1/p'}}.
$$
Therefore the function
$$
y_0(x)= e^{\;\int\limits_a^x \frac{dt}{w(t)}}\quad (a\;\;\mbox{be\; fixed \;in}\; (0, \infty))
$$
satisfies the problem (2.7), (2.8).

This completes the proof of Lemma 2.

Thus, we have the following
\begin{theorem} Let $1< p\le q(x)\le \overline q< \infty$ and $K< +\infty.$ Suppose that $\omega_1$
and $\omega_2$ are weight functions defined on $(0, \infty)$ and the derivative
$\omega'_1(t)$ exists for all $t\in (0, \infty)$ and $\omega_1(t)\ge \omega_1(0)> 0.$
Then the following statements are equivalent:

a)$\;$ there is a positive solution of the equation
$$
L\left(t,\omega_2,y\right)- \lambda\omega_1(t) \left(y'(t)\right)^{1/p'}= 0,
$$
$$
y(t)> 0,\,\;y'(t)> 0,\;\, y\in AC(0,\,\infty),\;
$$
where $\lambda> 0;$

b)$\;$ the weighted norm inequality
$$
\|u\|_{L_{q(\cdot), \omega_2}(0, \infty)}\le C_0\|u'\|_{L_{p, \omega_1}(0, \infty)} \eqno(2.16)
$$
holds,  where $u\in AC(0,\,\infty),$ $u(0)= \lim\limits_{t\to +0} u(t)= 0$ and $C_0> 0$ is independent of $u.$
\end{theorem}

\begin{remark}
Note that for $q(x)= q= const$ and $1< p\le q< \infty$ Theorem 4 was proved in [14]. Further, in [6], [44] and [45]
for $q(x)= q= const$ the one-dimensional two-weighted Hardy inequality and its equivalent differential
form (2.16) was connected with the Euler-Lagrange differential equation. Note that the proof of Theorem 4 is
based on the paper [14].
\end{remark}

{\bf Acknowledgement.} This work was supported by the Science Development Foundation under the President of
the Republic of Azerbaijan EIF-2010-1(1)-40/06-1.

\newpage
\begin{center}
{\bf References}
\end{center}

[1] K.F.Andersen, B.Muckenhoupt,  Weighted weak type Hardy inequalities with applica-tions to Hilbert transform
and maximal function, Studia Math., (1){\bf 72}(1982)  9-26.

[2] J.Appell, A.Kufner,  On the two-dimensional Hardy operator in Lebesgue spaces with mixed norms. Analysis,
15(1995) 91-98.

[3] R.A.Bandaliev,  On an inequality in Lebesgue space with
mixed norm and with variable summability exponent, Mat. Zametki, {\bf 3}
(84)(2008), 323-333.(In Russian). English translation:  Math. Notes,
{\bf 3}(84)(2008), 303-313 (2008).

[4] R.A.Bandaliev,  The boundedness of certain sublinear operator in the weighted variable
Lebesgue spaces, Czechoslovak Math. J., (2){(\bf 60)}(2010)  327-337.

[5] R.A.Bandaliev, The boundedness of multidimensional Hardy operator in the weighted variable
Lebesgue spaces, Lithuanian Math. J., (3){(\bf 50)}(2010) 249-259.

[6] P.R.Beesack,  Hardy's inequality and its extensions, Pacific J. Math., {\bf 11}(1961) 39-61.

[7] J.Bradley,  Hardy inequalities with mixed norms, Canadian Mathematical Bull. {\it 21}(1978) 405-408.

[8] A. \v{C}i\v{z}me\v{s}ija, J. Pe\v{c}ari\'{c}, I. Peri\'{c}, Mixed means and inequalities of Hardy
and Levin-Cochran-Lee type for multidimensional balls, Proc. Amer. Math. Soc., 9({\bf 128})(2000) 2543-2552.

[9] \v{C}i\v{z}me\v{s}ija A., Persson L.-E. and Wedestig A.,  Weighted integral inequalities of Hardy and
geometric mean operators with kernels over cones in $R^n,$  Ital. J. Pure Appl. Math., {\bf 18}(2005) 89-118.

[10] J.A.Cochran and C.-S.Lee,  Inequalities related to Hardy's and Heinig's, Math. Proc. Camb. Phil.Soc.,
{\bf 96}(1984) 1-7.

[11] L.Diening, P.Harjulehto, P.H\"{a}st\"{o},   and  M. R\.{u}\v{z}i\v{c}ka,
 Lebesgue and Sobolev spaces with variable exponents, Springer Lecture Notes, v.2017,
Springer-Verlag, Berlin, 2011.

[12] P. Dr\'{a}bek, H.P.Heinig  and A. Kufner,  Higher dimensional Hardy inequality,
Inter. Ser. Num. Math., {\bf 123}(1997)  1-16.

[13] B. Gupta, P.Jain, L.-E.Persson, A. Wedestig,  Weighted geometric mean inequalities over cones in $R^n,$
J. Inequal. Pure Appl. Math., art. 68(4)({\bf 4})(2003).

[14] P. Gurka,   Generalized Hardy's inequality, \v{C}as. pro P\^{e}stov\'{a}n\'{\i} Mat.{\bf 109}(1984)
 194-203.

[15] G.H.Hardy, J.E. Littlewood,  G.P\'{o}lya,  Inequalities, Cambridge University Press, 1934.

[16] P.Harjulehto, P.H\"{a}st\"{o} and M.Koskenoja,  Hardy's inequality in a variable exponent
Sobolev space,  Georgian Math.J. {\bf 12}(3),  431-442 (2005).

[17] H.P.Heinig,  Some extensions of inequalities, SIAM J. Math. Anal., ({\bf 6})(1975) 698-713.

[18] H.P. Heinig, R.Kerman, M.Krbec,  Weighted exponential inequalities, Georgian Math. J.,
8({\bf 1})(2001) 69-85.

[19] P.Jain, R.Hassija,  Some remarks on two dimensional Knopp type inequalities, Appl. Math. Letters,
 16({\bf 4})(2003) 459-465.

[20] P.Jain, L.E.Persson, A.Wedestig,  From Hardy to Carleman and general mean-type inequalities,
Function Spaces and Appl., 2000  117-130.

[21] P.Jain, L.E.Persson, A.Wedestig,  Carleman-Knopp type inequalities via Hardy's inequality,
Math. Inequal. Appl. (3){\bf 4}(2001) 343-355.

[22] K.Knopp,  \"{U}ber reihen mit positivern gliedern, J. London Math. Soc., {\bf 3}(1928) 205-211.

[23] V.M. Kokilashvili,  On Hardy's inequality in weighted spaces, Bull. Acad. Sci. Georgian SSR,
(1)({\bf 96})(1978) 37-40 (in Russian).

[24] T.S. Kopaliani,  On some structural properties of Banach function spaces and boundedness of certain
integral operators, Czechoslovak Math. J. (54){\bf 129}(2004) 791-805.

[25] O. Kov\'{a}\v{c}ik, J. R\'{a}kosn\'{\i}k,  On spaces $L^{p(x)}$ and $W^{k, p(x)},$
Czechoslovak Math. J. (41){\bf 116}(1991) 592-618.

[26] A.Kufner, H.Triebel, Generalization of Hardy's inequality, Conf. Semin. Mat.Univ. Bari
{\bf 156}(1978) 1-21.

[27] E.R. Love, Inequalities related to those of Hardy and of Cochran and Lee, Math. Proc. Camb. Phil.Soc.,
{\bf 99}(1986) 395-408.

[28] F.I.Mamedov, A.Harman,  On a weighted inequality of Hardy type in spaces $L^{p(\cdot)},$
J. Math. Anal. Appl., {\bf 353}(2009) 521-530.

[29] M.Marcus, V. Mizel, Y. Pinchover,  On the best constant for Hardy's inequality in $R^n,$
Trans. Amer. Math. Soc., 8({\bf 350})(1998), 3237-3255.

[30] V.G.Maz'ya,  Sobolev spaces,  Springer-Verlag, Berlin-Heidelberg-New York, 1985.

[31] B. Muckenhoupt,  Hardy's inequalities with weights in two dimensions, Notices Amer. Math. Soc., {\bf 29}(1982) 538.

[32]  J. Musielak,  Orlicz spaces and modular spaces, Lecture Notes in Math. 1034. Springer-Verlag,
Berlin-Heidelberg-New York, 1983.

[33] H.Nakano,  Modulared semi-ordered linear spaces, Maruzen, Co., Ltd., Tokyo 1950.

[34] B.Opic, P.Gurka, $N$-dimensional Hardy's inequality and embedding theorems for
weighted Sobolev spaces on unbounded domains,  Function spaces, Differential operator and Nonlinear
analysis (Sodankyl\"{a}, 1988) 108-124, Pitman Res. Notes Math. Ser., {\bf 211}, Longman Sci. Tech.,
Harlow, 1989.

[35] B. Opic, P.Gurka,  Weighted inequalities for geometric means,
Proc. Amer. Math. Soc., (3){\bf 120}(1994) 771-779.

[36] W. Orlicz,   \"{U}ber konjugierte exponentenfolgen,  Studia Math.{\it 3}(1931) 200-212.

[37] L.-E. Persson, V.D. Stepanov,  Weighted integral inequalities with the
geometric mean operator, J. Inequal. Appl., (5){\bf 7}(2002) 727-746.

[38] Pick L., Opic B.,  On geometric mean operator, J. Math. Anal. Appl., (3){\bf 183}(1994) 652-662.

[39] K.R.Rajagopal, M. R\.{u}\v{z}i\v{c}ka,  Mathematical modeling of  electrorheological materials, Cont. Mech.
and Termodyn., {\bf 13}(2001) 59-78.

[40] Y. Rakotondratsimba,  Weighted inequalities for the two-dimensional Laplace transform, SUT
 J. Math., (1){\bf 35}(1999) 11-36.

[41] E.T. Sawyer,  Weighted inequalities for the two-dimensional Hardy operator, Studia Math., ({\bf 82})(1985)
1-16.

[42] I.I. Sharapudinov,  On a topology of the space $L^{p(t)}([0,1]),$  Math. Notes {\it 26}(1979) 613-637.

[43] G. Sinnamon, One-dimensional Hardy type inequalities in many dimensions, Proc. Roy. Soc. Edinburgh,
sect A(4){\bf 128}(1998) 833-848.

[44] G. Talenti,  Osservazione sopra una classe di disuguaglianze, Rend. Sem. Mat. Fiz. Milano, {\bf 39}(1969)
171-185.

[45] G. Tomaselli,  A class of inequalities,  Boll.Un.Mat.Ital., {\bf 2}(1969) 622-631.

[46] A. Wedestig,  Some new Hardy type inequalities and their limiting inequalities,  J.of Ineq. in Pure
and  Appl. Math., art. 61(3) {\bf 4} (2003).

[47] A. Wedestig,  Weighted inequalities for the Sawyer two-dimensional Hardy operator and its limiting geometric
mean operator,  J. Inequal. Appl., {\bf 4} (2005) 387-394.

[48] X. Xu, Y. An, Existence and multiplicity of solutions for elliptic systems with nonstandard growth condition,
Nonlinear Analysis TMA, (4){\bf 68}(2008) 956-968.

[49] Q.H.Zhang, Existence and asymptotic behavior of positive solutions for variable exponent elliptic systems,
Nonlinear Analysis TMA, (1){\bf 70}(2009) 305-316.

[50] V.V. Zhikov,  Averaging of functionals of the calculus of variations and elasticity theory. Izv.  Akad. Nauk SSSR.
{\bf 50}(1986) 675-710. (In Russian). English transl.: Math. USSR, Izv., {\bf 29}(1987) 33-66.

\end{document}